\documentclass[11pt]{amsart}

\usepackage[T1]{fontenc}
\usepackage[utf8]{inputenc}
\usepackage{amsmath,amssymb,amsthm,mathtools}
\usepackage{geometry}
\usepackage{hyperref}
\usepackage{enumitem}
\usepackage{color}
\newtheorem{theorem}{Theorem}[section]
\newtheorem{proposition}[theorem]{Proposition}
\newtheorem{lemma}[theorem]{Lemma}
\newtheorem{corollary}[theorem]{Corollary}
\theoremstyle{remark}
\newtheorem{remark}[theorem]{Remark}

\newcommand{\Hessp}{\operatorname{Hess}^{+}}
\newcommand{\R}{\mathbb{R}}
\newcommand{\ip}[2]{\left\langle #1,#2\right\rangle}
\newcommand{\norm}[1]{\left\lVert #1\right\rVert}

\title{High-level convexity for products of squared Euclidean distance functions}
\author{Tudor Micu}
\address{Babe\c{s}-Bolyai University, Cluj-Napoca, Romania}
\email{tudor.micu@ubbcluj.ro}

\author{Cornel Pintea}
\address{Babe\c{s}-Bolyai University, Cluj-Napoca, Romania}
\email{cornel.pintea@ubbcluj.ro}

\author{George C. Turca\c{s}}
\address{Babe\c{s}-Bolyai University, Cluj-Napoca, Romania}
\email{george.turcas@ubbcluj.ro}

\date{\today}

\begin{document}

\begin{abstract}
We study smooth functions on Euclidean space whose Hessian is positive definite outside a bounded set, with emphasis on products of squared distance functions. More precisely, we first prove a simple convexity principle: if the superlevel region $f^{-1}([c,\infty))$ is contained in the Hessian-positive region of $f$, then the sublevel set $\{f\le c\}$ is convex. We apply this to finite products $F_P(x)=\prod_{p\in P}\norm{x-p}^2$, proving that their Hessian-positive complements are bounded. For the two-centre product $F_{p,q}(x)=\norm{x-p}^2\norm{x-q}^2$ in dimension $n\ge2$, we compute the Hessian-positive region and the exact value
\[
 h_{\max}(F_{p,q})=\frac{\norm{p-q}^4}{4}.
\]
This value is sharp for convexity of sublevel sets in the following sense: we prove convexity above it and nonconvexity below it. This also gives the exact convexity and quasiconvexity truncation levels for the two-centre model.
\end{abstract}

\maketitle

\section{Introduction}

The set where the Hessian of a smooth function is positive definite could be regarded as a tool for discussing the convexity of the function near infinity. From a broader convex-analysis viewpoint, this also touches the Fenchel problem of level sets, namely the problem of recognising convexity from the geometry of level families \cite{Rapcsak2005FenchelLevelSets}. In this short note, we
use this tool to study sublevel sets, continuing the line of high-level-set and
product results in
\cite{Pintea2022LevelSets,BrojbeanuPintea-JCA,BrojbeanuPintea2022HessPlus}, and
apply it to products of squared Euclidean distances. The main class of examples is
\[
 F_P(x)=\prod_{i=1}^m \norm{x-p_i}^2,
\]
where $p_i\in\R^n$, repetitions among the centres are allowed and $m\ge1$.

The related literature studies functions with bounded Hessian-positive
complement and their high level sets, as well as products that preserve bounded
Hessian-positive complement. In particular, Brojbeanu and Pintea studied
products with bounded $\Hessp$ complements
\cite{BrojbeanuPintea-JCA,BrojbeanuPintea2022HessPlus}. The present article
isolates the finite products of squared Euclidean distances and proves the
boundedness of their Hessian-positive complements directly, for all finite
configurations, by a leading-term argument. The final truncation-level
statement uses Pintea's truncation framework \cite{Pintea2026Cqt}; the
finite-product theorem and the two-centre Hessian computation are proved
directly in the present work.

The note has three main contributions. First, Proposition
\ref{prop:convexity-criterion} gives a direct dimension-free convexity
principle: if the high superlevel region of a $C^2$ function is contained in
its Hessian-positive region, then the corresponding sublevel set is convex. 
This isolates the fixed-level convexity mechanism behind earlier high-level
results. In contrast with the classical global first-/second-order criteria for quasi- and pseudo-convexity of $C^2$ functions, Proposition \ref{prop:convexity-criterion} is a level-wise sufficient criterion that only requires Hessian positivity on the relevant high-value region \cite{CrouzeixFerland1982Criteria}. Brojbeanu and Pintea prove such a high-level theorem for coercive
functions with bounded critical set and bounded Hessian-positive complement,
with convexity in the planar case \cite[Theorem 3.7]{BrojbeanuPintea2022HessPlus}.
Pintea proves compact connected regular high levels in all dimensions
\cite[Theorem~2.1]{Pintea2022LevelSets}, and the three-dimensional
ovaloid/convex-boundary conclusion under bounded Hessian-positive complement
\cite[Theorem~3.4 and Corollary~3.6]{Pintea2022LevelSets}. Our criterion does
not use planar curvature or three-dimensional ovaloid theory; it applies in
any dimension at every level satisfying
$f^{-1}([c,\infty))\subset\Hessp(f)$. Proposition
\ref{prop:bounded-critical-set} and Corollary
\ref{cor:high-level-consequences} record the resulting high-level
consequences.

Second, Theorem \ref{thm:finite-products-bounded-complement} gives the direct
leading-term proof that every nonempty finite product of squared distances has
bounded Hessian-positive complement. The corresponding high-level consequences
for this class are stated in Corollary \ref{cor:finite-products-high-levels}.
Third, and as the main model theorem, we compute the two-centre product
exactly.

For this model, writing
$c=(p+q)/2$, $u=(p-q)/2$, $a=\norm{u}$, and $y=x-c$, we compute the
normal form, gradient, Hessian, critical set, Hessian-positive region, and
compact Hessian-positive complement; see Lemma \ref{lem:two-centre-hessian},
Proposition \ref{prop:two-centre-critical-set}, Theorem
\ref{thm:two-centre-hess-region}, and Corollary
\ref{cor:two-centre-complement-bound}. The boundary identity in Theorem
\ref{thm:two-centre-hmax} gives, for every $n \geq 2$,
\[
 h_{\max}(F_{p,q})=4a^4=\frac{\norm{p-q}^4}{4}.
\]
The same value is the sharp infimum level for convexity of sublevel sets: we
prove convexity for $\lambda>4a^4$ in Corollary
\ref{cor:two-centre-high-sublevels-convex} and nonconvexity for
$0\le\lambda<4a^4$ in Theorem \ref{thm:two-centre-sharpness}. No endpoint
assertion for the sublevel $\{F_{p,q}\le4a^4\}$ is needed. In the sense of
\cite{Pintea2026Cqt}, Corollary \ref{cor:two-centre-truncation-levels} gives
the exact quasiconvexity and convexity truncation levels. For arbitrary finite
configurations we prove that $h_{\max}$ is finite and attained whenever the
Hessian-positive complement is nonempty, but we do not claim a closed formula;
that remains a quantitative problem.

\section{Convexity from the Hessian-positive region}

For a $C^2$ function $f:\R^n\to\R$, write
\[
 \Hessp(f)=\{x\in\R^n : Hf(x) \text{ is positive definite}\}.
\]
We say that $f$ is norm-coercive if $f(x)\to+\infty$ as $\norm{x}\to\infty$.
When
\[
 K_f:=\R^n\setminus\Hessp(f)
\]
is bounded, it is compact because it is closed, and we put
\[
 h_{\max}(f)=\sup_{x\in K_f} f(x),
\]
with the convention $\sup\emptyset=-\infty$. Thus, if $K_f$ is nonempty,
$h_{\max}(f)$ is a maximum.

The first observation we present is a convexity criterion. The proof does not require coercivity or any information about critical points. We just observe that a segment leaving a sublevel set would force an interior maximum along that segment, and such a maximum is incompatible with positive definiteness of the Hessian.

\begin{proposition}[Convexity criterion]\label{prop:convexity-criterion}
Let $f:\R^n\to\R$ be $C^2$. If
\[
 f^{-1}([c,\infty))\subseteq \Hessp(f),
\]
then the sublevel set $E_c=\{x:f(x)\le c\}$ is convex.
\end{proposition}

\begin{proof}
Let $a,b\in E_c$ and suppose the segment $\gamma(t)=(1-t)a+tb$ leaves $E_c$.
Then $\varphi=f\circ\gamma$ has an interior maximum $t_0$ with
$\varphi(t_0)>c$. Hence $\varphi''(t_0)\le0$. But
$\gamma(t_0)\in\Hessp(f)$, so
\[
 \varphi''(t_0)=(b-a)^T Hf(\gamma(t_0))(b-a)>0,
\]
contradiction.
\end{proof}

The above criterion separates the convexity from the usual compactness and regularity issues. To turn it into a high-level statement, we also need to know that critical points cannot occur arbitrarily \textit{far out}. This is already forced by the boundedness of the Hessian-positive complement, as we now show. 

\begin{proposition}[Bounded critical set]\label{prop:bounded-critical-set}
If $\R^n\setminus\Hessp(f)$ is bounded, then the critical set $C(f)=\{x:\nabla f(x)=0\}$ is bounded.
\end{proposition}

\begin{proof}
Choose $R\ge0$ such that $
 \R^n\setminus\Hessp(f)\subset B(0,R)$.

Suppose, toward a contradiction, that $C(f)$ is unbounded. Then there exists a sequence $(x_k)$ of critical points with $\norm{x_k}\to\infty$. Passing to a subsequence, we
may assume that $r_k=\norm{x_k}$ is strictly increasing and further that
$u_k=x_k/r_k$ converges in the unit sphere.

Fix $0<\delta<1$. Choose $k$ large enough so that
\[
 r_k\sqrt{1-\delta/2}>R.
\]
Then choose $l>k$ so large that
\[
 \ip{u_k}{u_l}\ge1-\delta.
\]
Put $x=x_k$, $y=x_l$, and $z(t)=(1-t)x+ty$. Since $r_l\ge r_k$, for
$0\le t\le1$ we have
\[
 \norm{z(t)}^2
 \ge r_k^2\bigl((1-t)^2+t^2+2t(1-t)(1-\delta)\bigr)
 =r_k^2\bigl(1-2\delta t(1-t)\bigr)
 \ge r_k^2(1-\delta/2)>R^2.
\]
Hence $z(t)\in\Hessp(f)$ for every $t\in[0,1]$.

Let $g(t)=f(z(t))$. Then
\[
 g''(t)=(y-x)^T Hf(z(t))(y-x)>0
\]
on $[0,1]$, so $g'$ is strictly increasing. But $x$ and $y$ are critical
points, and therefore
\[
 g'(0)=\ip{\nabla f(x)}{y-x}=0
 \text{ and }
 g'(1)=\ip{\nabla f(y)}{y-x}=0,
\]
a contradiction.
\end{proof}

We now combine the two elementary results above with coercivity. Above the height of the Hessian-positive complement, Proposition \ref{prop:convexity-criterion} gives convexity; above the critical values, the corresponding levels are regular. 

\begin{corollary}[High-level consequences]\label{cor:high-level-consequences}
Let $f:\R^n\to\R$ be $C^2$ and norm-coercive, and suppose that
$K_f=\R^n\setminus\Hessp(f)$ is bounded. Set
\[
 \gamma(f)=\max_{x\in C(f)} f(x).
\]
Then, for every
\[
 c>\max\{h_{\max}(f),\gamma(f)\},
\]
the sublevel set $E_c=\{f\le c\}$ is compact and convex, and the level set
$M_c=f^{-1}(c)$ is a compact regular $C^2$ hypersurface. If $n\ge2$, then
$M_c$ is connected and equals the boundary of a compact convex body.
\end{corollary}

\begin{proof}
By Proposition \ref{prop:bounded-critical-set}, $C(f)$ is bounded. It is closed,
and since $f$ is norm-coercive it is nonempty because $f$ attains a global
minimum. Hence $\gamma(f)$ is finite.

The set $E_c$ is closed and bounded, hence compact. Since $c>h_{\max}(f)$, every
point with $f(x)\ge c$ lies in $\Hessp(f)$. Proposition
\ref{prop:convexity-criterion} therefore implies that $E_c$ is convex.

Since $c>\gamma(f)$, the level set $M_c$ contains no critical point. Thus $c$ is a
regular value and $M_c$ is a $C^2$ embedded hypersurface. The global minimum
value of $f$ is strictly below $c$, and coercivity gives points with value above
$c$ along every sufficiently long ray; hence $M_c$ is nonempty and compact.
Moreover $M_c=\partial E_c$: points with $f<c$ are interior to $E_c$, and at a
regular point of the level $f=c$ there are nearby points on both sides of the
level.

It remains only to record connectedness when $n\ge2$. Let
$x_0\in\operatorname{int}(E_c)$. For $v\in S^{n-1}$, put
\[
 \rho(v)=\max\{t\ge0:x_0+tv\in E_c\}.
\]
Compactness gives existence of the maximum, and convexity shows that
$x_0+\rho(v)v$ is the unique point where the ray $x_0+\R_{\ge0}v$ meets
$\partial E_c$. The map $v\mapsto \rho(v)$ is continuous: if $v_j\to v$ and
$\rho(v_j)\to L$ along a subsequence, then $x_0+Lv\in E_c$, so
$L\le\rho(v)$; conversely, for every $t<\rho(v)$ the point $x_0+tv$ lies in
$\operatorname{int}(E_c)$, hence $x_0+tv_j\in E_c$ for all large $j$, giving
$\liminf_j\rho(v_j)\ge t$ and then $\liminf_j\rho(v_j)\ge\rho(v)$. Thus the
radial map $v\mapsto x_0+\rho(v)v$ is continuous and onto $\partial E_c$.
Since $S^{n-1}$ is connected for $n\ge2$, $\partial E_c$ is connected.
\end{proof}

\section{Finite products of squared distances}

Let $n\ge1$ and let $P=(p_1,\dots,p_m)$ be a nonempty finite list of points in
$\R^n$. Repeated centres are allowed. Set
\[
 F_P(x)=\prod_{i=1}^m \norm{x-p_i}^2.
\]

We first show that these have bounded Hessian-positive complements. The proof is purely asymptotic: \textit{at infinity}, $F_P$ has a certain leading term whose Hessian dominates all lower degree contributions. 

\begin{theorem}\label{thm:finite-products-bounded-complement}
For every nonempty finite list $P$ as above,
\[
 \R^n\setminus\Hessp(F_P)
\]
is bounded.
\end{theorem}

\begin{proof}
If $m=1$, then $HF_P=2I$ everywhere, so the complement is empty. Assume
$m\ge2$.

Write
\[
 F_P(x)=\norm{x}^{2m}+Q(x),\text{ where } \deg Q\le 2m-1.
\]
The Hessian of the leading term is
\[
 H\norm{x}^{2m}=2m\norm{x}^{2m-2}I+4m(m-1)\norm{x}^{2m-4}xx^T.
\]
Thus for every unit vector $v$,
\[
 v^T H\norm{x}^{2m}v\ge 2m\norm{x}^{2m-2}.
\]
Since every entry of $HQ$ has degree at most $2m-3$, there is a constant
$C>0$ such that
\[
 \norm{HQ(x)}_{\mathrm{op}}\le C(1+\norm{x})^{2m-3}
\]
for all $x$. Hence, uniformly over all unit vectors $v$,
\[
 v^T HF_P(x)v
 \ge 2m\norm{x}^{2m-2}-C(1+\norm{x})^{2m-3}.
\]
This lower bound is positive for $\norm{x}$ sufficiently large. Thus $HF_P$
is positive definite outside a ball, so the complement of $\Hessp(F_P)$ is
bounded.
\end{proof}

The argument above shows boundedness but does not yet reflect the geometry of the configuration. Recentering at the centroid removes the degree $2m-1$ term and makes the first configuration-dependent correction visible. 

\begin{proposition}[Centroid expansion]\label{prop:centroid-expansion}
Let
\[
 \bar p=\frac1m\sum_{i=1}^m p_i,\text{ and } p_i=\bar p+u_i,
\]
so that $\sum_i u_i=0$. Then 
\[
 F_P(\bar p+y)=\norm{y}^{2m}+\widetilde Q(y),
 \text{ where } \deg \widetilde Q\le 2m-2.
\]
If $m\ge2$, then the homogeneous term of degree $2m-2$ is governed by the
second-moment matrix $
 M_P=\sum_{i=1}^m u_i u_i^T$, More precisely, if we denote its trace by $
 \sigma_P=\operatorname{tr}(M_P)=\sum_{i=1}^m\norm{u_i}^2$, then 
\[
 F_P(\bar p+y)
 =
 \norm{y}^{2m}
 +\norm{y}^{2m-4}
  \bigl(\sigma_P\norm{y}^2-2y^T M_Py\bigr)
 +R(y),
\]
where $\deg R\le 2m-3$.
\end{proposition}

\begin{proof}
For $x=\bar p+y$ we have
\[
 \norm{x-p_i}^2=\norm{y-u_i}^2
 =\norm{y}^2-2\ip{y}{u_i}+\norm{u_i}^2.
\]
Thus
\[
 F_P(\bar p+y)
 =\prod_{i=1}^m
  \bigl(\norm{y}^2-2\ip{y}{u_i}+\norm{u_i}^2\bigr).
\]
The degree $2m$ term is $\norm{y}^{2m}$. The homogeneous term of degree
$2m-1$ comes from choosing the linear term in exactly one factor and
$\norm{y}^2$ in all the others; it is
\[
 -2\norm{y}^{2m-2}\sum_{i=1}^m\ip{y}{u_i}
 =-2\norm{y}^{2m-2}\ip{y}{\sum_i u_i}=0.
\]
This proves the first assertion.

Assume now that $m\ge2$. The homogeneous terms of degree $2m-2$ have two
sources. Choosing the constant term in one factor gives
$\sigma_P\norm{y}^{2m-2}$. Choosing the linear terms in two distinct factors
gives
\[
 4\norm{y}^{2m-4}\sum_{i<j}\ip{y}{u_i}\ip{y}{u_j}.
\]
Since $\sum_i\ip{y}{u_i}=0$, we have
\[
 \sum_{i<j}\ip{y}{u_i}\ip{y}{u_j}
 =-\frac12\sum_i\ip{y}{u_i}^2
 =-\frac12 y^T M_Py.
\]
Therefore the full homogeneous term of degree $2m-2$ is
\[
 \norm{y}^{2m-4}
 \bigl(\sigma_P\norm{y}^2-2y^T M_Py\bigr).
\]
All remaining terms have degree at most $2m-3$.
\end{proof}

We can now return from the asymptotic computations to the level-set geometry. Since $F_P$ is coercive and its Hessian-positive complement is bounded, the general results of the previous section apply directly.

\begin{corollary}\label{cor:finite-products-high-levels}
For every nonempty finite list $P$ in $\R^n$, the critical set $C(F_P)$ is
bounded and $h_{\max}(F_P)$ is finite and attained if
$\R^n\setminus\Hessp(F_P)$ is nonempty. Moreover, all sufficiently high
sublevel sets of $F_P$ are compact and convex. For $n\ge2$, all sufficiently
high levels are compact connected regular hypersurfaces bounding compact convex
bodies.
\end{corollary}

\begin{proof}
The boundedness of the critical set follows from Theorem
\ref{thm:finite-products-bounded-complement} and Proposition
\ref{prop:bounded-critical-set}. The same theorem also shows that the
Hessian-positive complement is compact if it is nonempty, so $h_{\max}(F_P)$ is
finite and attained.

Finally, $F_P$ is norm-coercive, since
$F_P(x)/\norm{x}^{2m}\to1$ as $\norm{x}\to\infty$. Applying Corollary
\ref{cor:high-level-consequences} gives the stated high-level conclusions.
\end{proof}

\section{The two-centre product}

Let $p\ne q$ and consider the function
\[
 F_{p,q}(x)=\norm{x-p}^2\norm{x-q}^2.
\]
Set
\[
 c=\frac{p+q}{2},\quad u=\frac{p-q}{2},\quad a=\norm{u},\quad y=x-c.
\]
Then
\[
 F_{p,q}(x)=(\norm{y}^2+a^2)^2-4\ip{y}{u}^2.
\]
Looking at this form, one can guess why the two-centre case is explicitly computable:
apart from the radial term, only the direction \(u\) plays a role.  The
quantity controlling the Hessian-positive region will be
\[
D(y)=3\|y\|^4+4\langle y,u\rangle^2-2a^2\|y\|^2-a^4.
\]

\begin{lemma}\label{lem:two-centre-hessian}
One has
\[
 \nabla F_{p,q}(x)=4(\norm{y}^2+a^2)y-8\ip{y}{u}u
\]
and
\[
 HF_{p,q}(x)=4(\norm{y}^2+a^2)I+8(yy^T-uu^T).
\]
\end{lemma}

\begin{proof}
Since $p=c+u$ and $q=c-u$, we have
\[
 F_{p,q}(x)=\norm{y-u}^2\norm{y+u}^2
 =\bigl(\norm{y}^2+a^2-2\ip{y}{u}\bigr)
  \bigl(\norm{y}^2+a^2+2\ip{y}{u}\bigr),
\]
which gives the displayed normal form. Differentiating with respect to $y$
gives the gradient formula, and a second differentiation gives the Hessian
formula. The translation $x\mapsto y=x-c$ does not change these derivatives.
\end{proof}

The gradient formula already determines the critical set.  Any component of
\(y\) orthogonal to \(u\) is multiplied by the positive factor
\(\|y\|^2+a^2\), so a critical point must lie on the line through the two
centres.

\begin{proposition}\label{prop:two-centre-critical-set}(see also \cite[Example 2.1]{BrojbeanuPintea-JCA})
The critical set of $F_{p,q}$ is
\[
 C(F_{p,q})=\{p,c,q\},
\]
and the corresponding critical values are $0,a^4,0$.
\end{proposition}

\begin{proof}
If $\nabla F_{p,q}(x)=0$, then
\[
 (\norm{y}^2+a^2)y=2\ip{y}{u}u.
\]
The component of $y$ orthogonal to $u$ must vanish, since
$\norm{y}^2+a^2>0$. Hence $y=tu$ for some $t\in\R$. Substitution gives
\[
 (t^2+1)t=2t,
\]
so $t\in\{-1,0,1\}$. Thus $x=c+tu$ is one of $p,c,q$, and direct evaluation
gives the stated critical values.
\end{proof}

For the computation of \(h_{\max}\), the critical points alone are not
enough; we need the full Hessian-positive region. The next result is a generalization to arbitrary dimension of \cite[Remark 3.2]{BrojbeanuPintea-JCA}. 

\begin{theorem}\label{thm:two-centre-hess-region}
For $n\ge2$,
\[
 HF_{p,q}(x)\succ0
 \quad\Longleftrightarrow\quad
 3\norm{y}^4+4\ip{y}{u}^2-2a^2\norm{y}^2>a^4.
\]
\end{theorem}

\begin{proof}
Let $e_1=u/a$. If $y$ is not parallel to $u$, choose a unit vector $e_2$
so that $y\in\operatorname{span}\{e_1,e_2\}$. If $y$ is parallel to $u$,
choose any unit vector $e_2\perp e_1$. Write $y=se_1+te_2$, with $t=0$
in the parallel case.

On $\operatorname{span}\{e_1,e_2\}^{\perp}$, the two rank-one terms in the
Hessian vanish, so $HF_{p,q}$ acts by the positive scalar
$4(\norm{y}^2+a^2)$. On $\operatorname{span}\{e_1,e_2\}$, the Hessian block is
\[
 4\begin{pmatrix}
 3s^2+t^2-a^2 & 2st\\
 2st & s^2+3t^2+a^2
 \end{pmatrix}.
\]
The determinant of the inner matrix is
\[
 3(s^2+t^2)^2+4a^2s^2-2a^2(s^2+t^2)-a^4
 =D(y),
\]
because $\ip{y}{u}=as$. The trace of the full $2\times2$ block is
$16(s^2+t^2)=16\norm{y}^2$. Thus, if $D(y)>0$, then $y\ne0$ and this trace is
positive; for a real symmetric $2\times2$ matrix this is equivalent to
positive definiteness of the block. Conversely, a positive definite block has
positive determinant, hence $D(y)>0$. Combining this with the positive
orthogonal directions proves the criterion. The same block decomposition gives
the useful determinant identity
\[
 \det HF_{p,q}(x)=4^n(\norm{y}^2+a^2)^{n-2}D(y),
\]
but the positive-definiteness criterion above uses the block decomposition, not
determinant positivity alone in higher dimension.
\end{proof}

The same inequality also gives a useful geometric bound.  Points in the
Hessian-positive complement cannot lie farther from the midpoint \(c\) than
the two centres themselves.

\begin{corollary}\label{cor:two-centre-complement-bound}
For $n\ge2$,
\[
 K:=\R^n\setminus\Hessp(F_{p,q})=\{x:D(x-c)\le0\}
\]
is compact and satisfies $K\subseteq\overline B(c,a)$.
\end{corollary}

\begin{proof}
If $x\in K$ and $r=\norm{y}$, then
\[
 3r^4+4\ip{y}{u}^2-2a^2r^2\le a^4.
\]
Dropping the nonnegative middle term gives
\[
 3r^4-2a^2r^2-a^4\le0.
\]
With $z=r^2$, this is
\[
 3z^2-2a^2z-a^4\le0,
\]
whose roots are $-a^2/3$ and $a^2$. Since $z\ge0$, we get $r\le a$. The set
$K=\{D\le0\}$ is closed, so it is compact.
\end{proof}

\section{Exact \texorpdfstring{$h_{\max}$}{hmax} and sharpness}

We now maximize \(F_{p,q}\) on the Hessian-positive complement.

\begin{theorem}\label{thm:two-centre-hmax}
For $p\ne q$ and $n\ge2$,
\[
 h_{\max}(F_{p,q})=4a^4=\frac{\norm{p-q}^4}{4}.
\]
\end{theorem}

\begin{proof}
Let $K=\R^n\setminus\Hessp(F_{p,q})$. By Corollary
\ref{cor:two-centre-complement-bound}, $K$ is compact, so $F_{p,q}|_K$
attains its maximum.

First suppose that a maximizer lies in the interior of $K$. Then it is an
unconstrained local maximum of $F_{p,q}$, hence a critical point. By
Proposition \ref{prop:two-centre-critical-set}, the only critical points are
$p,c,q$. The points $p$ and $q$ lie in $\Hessp(F_{p,q})$, since
$D(\pm u)=4a^4>0$, while $c\in K$, since $D(0)=-a^4$, and
$F_{p,q}(c)=a^4$. Hence every interior critical maximum on $K$ has value at
most $a^4$.

Now suppose that a maximizer lies on $\partial K$. Then $D(y)=0$, i.e.
\[
 3\norm{y}^4+4\ip{y}{u}^2-2a^2\norm{y}^2=a^4.
\]
Substitution into the normal form gives the boundary identity
\[
 F_{p,q}(x)=4\norm{y}^4.
\]
Since $K\subseteq\overline B(c,a)$, every boundary value is at most $4a^4$.
Equality is attained whenever $\norm{y}=a$ and $\ip{y}{u}=0$, which is possible
because $n\ge2$. Therefore
\[
 h_{\max}(F_{p,q})=4a^4=\frac{\norm{p-q}^4}{4}.
\]
\end{proof}

Once \(h_{\max}\) is known, convexity above this value follows directly from
the criterion of Section 2.

\begin{corollary}\label{cor:two-centre-high-sublevels-convex}
For every $\lambda>4a^4$, the sublevel set $\{F_{p,q}\le\lambda\}$ is convex.
\end{corollary}

\begin{proof}
If $F_{p,q}(x)\ge\lambda>h_{\max}(F_{p,q})$, then $x\in\Hessp(F_{p,q})$.
Thus $F_{p,q}^{-1}([\lambda,\infty))\subseteq\Hessp(F_{p,q})$, and
Proposition \ref{prop:convexity-criterion} applies.
\end{proof}

It remains to show that the constant cannot be lowered.  The construction
below exhibits, for every level below \(4a^4\), two points in the sublevel
set whose midpoint lies outside it.

\begin{theorem}[Sharpness below the threshold]\label{thm:two-centre-sharpness}
For $n\ge2$ and every $0\le\lambda<4a^4$, the sublevel set
$\{F_{p,q}\le\lambda\}$ is not convex.
\end{theorem}

\begin{proof}
Put $\rho^2=\lambda/(4a^2)$. Since $\lambda<4a^4$, we have $0\le\rho<a$.
Choose $w\perp u$ with $\norm{w}=\rho$, and set $e=u/a$. Define
\[
 y_\pm=\pm\sqrt{a^2-\rho^2}\,e+w.
\]
Then $\norm{y_\pm}=a$ and
\[
 \ip{y_\pm}{u}=\pm a\sqrt{a^2-\rho^2}.
\]
Using the normal form,
\[
 F_{p,q}(c+y_\pm)
 =(2a^2)^2-4a^2(a^2-\rho^2)
 =4a^2\rho^2
 =\lambda.
\]
Thus $c+y_+$ and $c+y_-$ lie in the sublevel set. Their midpoint is $c+w$,
and $\ip{w}{u}=0$, so
\[
 F_{p,q}(c+w)=(a^2+\rho^2)^2.
\]
Therefore
\[
 F_{p,q}(c+w)-\lambda
 =\left(a^2+\frac{\lambda}{4a^2}\right)^2-\lambda
 =\frac{(4a^4-\lambda)^2}{16a^4}>0.
\]
The midpoint is not in the sublevel set. Hence the sublevel set is not convex.
For $\lambda=0$, this is the same construction with $w=0$ and endpoints
$p,q$.
\end{proof}

\begin{remark}
Corollary \ref{cor:two-centre-high-sublevels-convex} and Theorem
\ref{thm:two-centre-sharpness} show that $4a^4$ is the sharp infimum level for
convexity of the sublevel sets. We do not use or assert convexity at the
endpoint $\lambda=4a^4$.
\end{remark}

Finally, we translate this sharp level statement into the truncation
terminology of \cite{Pintea2026Cqt}.  For a real number \(\tau\), write
\(T_\tau(f)=\max\{\tau,f\}\), and denote by \(sql(f)\) and \(scl(f)\) the
smallest quasiconvexity and convexity truncation levels, respectively.. The truncation levels are taken over
$\tau\ge\inf f$, as defined in \cite[Section 2]{Pintea2026Cqt}.

\begin{corollary}\label{cor:two-centre-truncation-levels}
For $n\ge2$,
\[
 \operatorname{sql}(F_{p,q})
 =\operatorname{scl}(F_{p,q})
 =h_{\max}(F_{p,q})
 =\frac{\norm{p-q}^4}{4}.
\]
\end{corollary}

\begin{proof}
Let $h=4a^4$. Corollary \ref{cor:two-centre-high-sublevels-convex} shows that
for every $\tau>h$, all sublevel sets $\{F_{p,q}\le r\}$ with $r\ge \tau$
are convex. The sublevel sets of $T_\tau(F_{p,q})$ below $\tau$ are empty.
Hence $T_\tau(F_{p,q})$ is quasiconvex for every $\tau>h$, and
$\operatorname{sql}(F_{p,q})\le h$.

Conversely, $\inf F_{p,q}=0$. Thus the only truncation levels below $h$ allowed
by the definition satisfy $0\le\tau<h$. For such a $\tau$, choose $r$ with
$\tau\le r<h$. By Theorem \ref{thm:two-centre-sharpness}, the sublevel set
$\{F_{p,q}\le r\}$ is not convex. Since $r\ge\tau$, this is also a sublevel set
of $T_\tau(F_{p,q})$, so $T_\tau(F_{p,q})$ is not quasiconvex. Thus
$\operatorname{sql}(F_{p,q})\ge h$, and
$\operatorname{sql}(F_{p,q})=h$.

The critical set of $F_{p,q}$ is finite by Proposition
\ref{prop:two-centre-critical-set}, and $\R^n\setminus\Hessp(F_{p,q})$ is
compact by Corollary \ref{cor:two-centre-complement-bound}. Therefore
\cite[Theorem 3.1]{Pintea2026Cqt} applies and gives
\[
 \operatorname{sql}(F_{p,q})\le \operatorname{scl}(F_{p,q})
 \le \max\{\operatorname{sql}(F_{p,q}),h_{\max}(F_{p,q})\}=h.
\]
Hence $\operatorname{scl}(F_{p,q})=h$ as well.
\end{proof}

\section{Cassini comparison and open problems}

For $p=(-a,0)$ and $q=(a,0)$ in the plane,
\[
 F_{p,q}(x,y)=((x+a)^2+y^2)((x-a)^2+y^2)=f_a(x,y)+a^4,
\]
where
\[
 f_a(x,y)=(x^2+y^2)^2-2a^2(x^2-y^2).
\]
Thus the classical Cassini polynomial $f_a$ differs from the squared-distance
product by the constant $a^4$. This shifted quartic is the planar example
used in \cite{PiTo,Pintea2026Cqt}. The shift matters for every threshold
statement: the product value at the midpoint is $a^4$, while the exact
Hessian-complement maximum for the product is $4a^4$.

The preceding finite-product theorem gives bounded Hessian-positive complement
for arbitrary nonempty finite configurations and hence guarantees that
$h_{\max}(F_P)$ is finite and attained whenever this complement is nonempty.
Determining it explicitly for general configurations with more than two centres
is not addressed here. Even in the plane, the problem becomes a polynomial
optimization problem on a compact set defined by polynomial inequalities;
special symmetric configurations may admit formulas, but no configuration-free
expression is asserted in this article.

\bibliographystyle{alpha}
\bibliography{references}
\end{document}